# SINGULARITIES OF MEAN CURVATURE FLOW OF SURFACES

# Preliminary Version

T. Ilmanen[1]

Department of Mathematics
University of Wisconsin
Madison, WI 53706

A family of $k$-dimensional surfaces $\{M_t\}_{t \in [0,T)}$ smoothly immersed in $\mathbf{R}^n$ is moving by mean curvature provided

(1) $$\frac{\partial}{\partial t} x = \vec{H}(x), \qquad x \in M_t, \quad t \in [0, T),$$

where $\vec{H} = (H^\alpha)_{\alpha=1}^n$ is the mean curvature vector of $M_t$ at $x$. (We will suppress explicit indication of the immersion map whenever possible.)

It is well known that the surface $M_t$ typically develops singularities at some finite time, which we take to be $T$. A major problem has been to study the nature of these singularities. Some classical examples include the shrinking sphere, the spirograph curves of Abresch-Langer [AL], the neckpinch of Grayson [Gr], the cusp singularity of Angenent [A1], the nipple of Angenent-Velazquez [AV2], and fully singular flows in Brakke's classic work [B].

---
[1]The author gratefully acknowledges the support of the NSF through a Postdoctoral Fellowship, Grant DMS-9205153, and the PYI Grant of S. Angenent at the University of Wisconsin.

The big breakthrough was Huisken's discovery [H2] of the now well-known *monotonicity formula* for motion by mean curvature, which states that

$$(2) \quad \frac{d}{dt} \int_{M_t} \rho(x,t) \, d\mathcal{H}^k(x)$$
$$= -\int_{M_t} \rho(x,t) \left| \vec{H}(x,t) + \frac{S(x,t)^\perp \cdot (x-y)}{2(s-t)} \right|^2 d\mathcal{H}^k(x),$$

for $t < s$, where $\rho = \rho_{y,s}$ is the backwards $k$-dimensional heat kernel defined by

$$\rho_{y,s}(x,t) = \frac{1}{(4\pi(s-t))^{k/2}} e^{-|x-y|^2/4(s-t)}, \quad t < s, \quad x, y \in \mathbf{R}^n,$$

$S(x,t)^\perp$ is the projection to the orthogonal complement of $T_x M_t$, and $\mathcal{H}^k$ denote Hausdorff measure on $M_t$ with respect to the induced metric. This formula is analogous to the monotonicity formula for minimal surfaces [F, 5.4.3], the monotonicity formula of Giga-Kohn [GK, Prop. 3], the mean value property for harmonic functions (see also the generalization by Chen-Struwe [CS, Lemma 2.2]), and for the Yang-Mills flow [P]. The formula leads us to rescale the flow parabolically about $(y,T)$ by defining, for $\lambda > 0$,

$$(3) \quad M_t^\lambda \equiv M_t^{(y,T),\lambda} = \lambda^{-1} \cdot (M_{T+\lambda^2 t} - y), \quad t \in [-T/\lambda^2, 0),$$

where $\lambda^{-1} \cdot$ represents homothetic expansion and $-y$ represents translation. Using (2), Huisken proves that, under the so-called *Type I* hypothesis (a natural curvature bound that implies a uniform, local $C^2$ bound for the rescaled surfaces), there exist subsequences $\{\lambda_i\}$ such that $\{M_t^{\lambda_i}\}_{t \in [-T/\lambda^2, 0)}$ converges locally in $\mathbf{R}^n \times (-\infty, 0)$, to a smooth limit $\{N_t\}_{t<0}$ that is self-similarly shrinking, that is,

$$(4) \quad N_t = \sqrt{-t} \cdot N_{-1}, \quad t < 0,$$

which (by (1) and separation of variables) is equivalent to

$$\vec{H}(x) + \frac{S(x)^\perp \cdot x}{2} = 0, \qquad x \in N_{-1}, \tag{5}$$

a parametric elliptic equation. The flow $\{N_t\}_{t<0}$ is called the *blowup* at $(y, T)$.

Known solutions of this equation include spheres, cylinders, the immersed curves of Abresch-Langer [AL], the shrinking donut of Angenent and various immersed surfaces of rotation [A2], minimal cones, new rotationally symmetric surfaces in $\mathbf{R}^4$ through $\mathbf{R}^7$ [AIV2], and new examples of Chopp [C] in $\mathbf{R}^3$.

It has been desirable to remove the Type I hypothesis, since it can be verified in only a few cases (see [H1, AAG] for example). We will show in complete generality (Lemma 8) that any family of rescalings $\{M_t^\lambda\}_{t \in [-T/\lambda^2, 0)}$ converges subsequentially to a self-similarly shrinking mean curvature flow $\{\nu_t\}_{t<0}$ in the weak sense of Brakke [B]. This result has also been obtained by White [W].

The question then becomes the regularity of this blowup. Let us assume for the duration that

$$\sup_{x \in \mathbf{R}^n} \sup_{R > 0} \frac{\mathcal{H}^k(M_0 \cap B_R(x))}{\omega_k R^k} \leq D, \tag{6}$$

where $\omega_k$ is the volume of the unit ball in $\mathbf{R}^k$. Our main result is the following.

**Theorem 1.** Let $\{M_t\}_{t \in [0,T)}$ be a family of properly immersed 2-dimensional surfaces flowing smoothly by mean curvature in $\mathbf{R}^n$. Suppose (6) holds and $M_0$ has finite genus.

Then for any $(y, T)$ and any weak blowup $\{\nu_t\}_{t<0}$ at $(y, T)$, spt $\nu_t$ is the image of a conformal branched immersion satisfying (5) in the parametrized sense. In particular, spt $\nu_t$ is the image of a smoothly immersed surface away from a discrete set $Q$ in $\mathbf{R}^n$.

Theorem 1 will be stated more precisely in §4. As an immediate corollary,

**Theorem 2.** *If $n = 3$ and $M_0$ is embedded, then spt $\nu_t$ is smooth.*

The reader should beware that the convergence to the blowup might not be smooth, even when the support of the blowup is smooth; there might, in principle, be small necks and other topological kinks that pinch off and disappear in the limit.

We would like to thank Sigurd Angenent and Klaus Ecker for very valuable discussions. Also, after proving Theorem 3 initially by a convoluted argument involving curve-shortening and tube-counting, the author consulted with R. Gulliver, R. Kusner, T. Toro, B. White, and R. Ye, whose suggestions led to the present form and proof of the Theorem.

**Main Idea.** The author was originally inspired by a talk given by K. Ecker, who proves that there is $\varepsilon_0 > 0$ such that if $\{M_t\}_{t \in [0,1)}$ is a flow of surfaces in $\mathbf{R}^3$ (or a 3-manifold) satisfying

$$\sup_{t \in [0,1)} \int_{M_t \cap B_1} |A|^2 \leq \varepsilon_0^2$$

then the flow can be extended smoothly up to $t = 1$, with the estimate

$$|A(0,1)| \leq C \int_0^1 \int_{M_t \cap B_1} |A|^2.$$

For a smooth, closed surface $M^2$ immersed in $\mathbf{R}^n$, the Gauss equation [Do, p.135] implies the relation

$$|\vec{H}|^2 - |\vec{A}|^2 = 2K,$$

where $\vec{A} = (A_{ij}^\alpha)$ is the second fundamental form, $\vec{H} = (H^\alpha) = (A_{ii}^\alpha)$ is the mean curvature vector, and $K = R_{ijij}^M/2$ is the Gauss curvature of $M$. Here $i, j = 1, 2$ and $\alpha = 1, \ldots, n$. Then the Gauss-Bonnet formula says

$$\int_M |\vec{A}|^2 = \int_M |\vec{H}|^2 - 2K = \int_M |\vec{H}|^2 - 4\pi\chi(M),$$

where $\chi(M)$ is the Euler characteristic. Since $\chi(M_t)$ is unchanging, the usual formula for area decrease [H1, Cor. 3.6] implies the curvature estimate

$$\int_0^T \int_{M_t} |\vec{A}|^2 \leq \mathcal{H}^2(M_0) - 4\pi T \chi(M_0).$$

Unfortunately, the estimate deteriorates in usefulness as we approach $T$, reflecting the fact that the left-hand side is not a scale-invariant quantity. In fact, in a neck-pinch, or any of the noncompact examples of Chopp [C], $\int_{M_t \cap U} |\vec{A}|^2 \to \infty$ as $t \to T$ in each neighborhood $U$ of the singularity.

By (2) and (6) we have $\iint \rho \leq D$ at all points and all scales, so for mean curvature flows of any dimension,

(7) $$\mathcal{H}^k(M_t \cap B_r(x)) \leq C(k)Dr^k, \qquad t \in [0, T), \quad x \in \mathbf{R}^n.$$

Then (2) implies the scale-invariant mean curvature estimate

(8) $$r^{-k} \int_{t-r^2}^t \int_{M \cap B_r} |\vec{H}|^2 \leq C(k)D$$

which applies equally well to the rescaled flows $\{M_t^\lambda\}_{t \in [-T/\lambda^2, 0]}$.

In order for this to be useful, we must localize the Gauss-Bonnet estimate. If a 2-manifold $M$ has a smooth boundary $\partial M$, then

(9) $$\begin{aligned} \int_M |\vec{A}|^2 &= \int_M |\vec{H}|^2 - 4\pi\chi(M) + \int_{\partial M} 2\tilde{k} \cdot n \\ &= \int_M |\vec{H}|^2 - 4\pi(2c(M) - 2g(M) - h(M)) + \int_{\partial M} 2\tilde{k} \cdot n, \end{aligned}$$

where $c(M)$ is the number of components of $M$, $g(M)$ is the genus of $M$, defined to be the genus of the closed surface obtained by capping off the boundary components of $M$ by disks, $h(M)$ is the number of components of $\partial M$, $\tilde{k}$ is the curvature vector of $\partial M$ in $M$, and $n$ is the inward pointing unit co-normal of $M$ along $\partial M$.

Replacing $M$ by $M \cap B_r$ and integrating over $r \in [1, 2]$, in §1 we control the two boundary terms to yield the following central estimate.

**Theorem 3 (Local Gauss-Bonnet Estimate).** Let $R > 1$ and let $M$ be a 2-manifold properly immersed in $B_R$. Then for any $\varepsilon > 0$,

$$(10) \quad (1-\varepsilon) \int_{M \cap B_1} |\vec{A}|^2 \, d\mathcal{H}^2 \leq \int_{M \cap B_R} |\vec{H}|^2 \, d\mathcal{H}^2 + 8\pi g(M \cap B_R)$$

$$- 8\pi c'(M \cap B_R) + \frac{24\pi D' R^2}{\varepsilon (R-1)^2}$$

where

$$D' = \sup_{S \in [1,R]} \frac{\mathcal{H}^2(M \cap B_S)}{\pi S^2}$$

and $c'(M \cap B_R)$ is the number of components of $M \cap B_R$ that meet $B_1$.

The above estimate should be useful in numerous geometric problems involving two-dimensional surfaces.

Integrating this estimate in time, together with (8), implies the following (see §1).

**Theorem 4 (Integral Curvature Estimate).** Let $\{M_t\}_{t \in [0,T)}$ be a mean curvature flow of 2-manifolds smoothly and properly immersed in $\mathbf{R}^n$. Then for every $B_r(x) \times [t - r^2, t] \subseteq \mathbf{R}^n \times [0, T)$ and $\varepsilon > 0$,

$$(11) \quad (1-\varepsilon) r^{-2} \int_{t-r^2}^{t} \int_{M_t \cap B_r(x)} |\vec{A}|^2 \leq 8\pi g(M_0) + \frac{CD}{\varepsilon}.$$

In §2 we make a few technical points needed to prove Theorem 1. The first, mentioned above, is that the blowup procedure of Huisken always works in a weak sense, yielding a varifold solution $\nu_{-1}$ of (5). Allard's Regularity Theorem implies that if $\nu_{-1}$ consists of roughly a single layer near a point $x$, then $\mathrm{spt}\,\nu_{-1}$ is smooth near $x$. A lemma of Simon implies that if $\int_{M_t^\lambda \cap B_r(x)} |\vec{A}|^2$ is sufficiently small, then $M_t^\lambda \cap B_{\sigma r}(x)$ decomposes as a union of nearly flat, embedded disks (layers) for near $x$. Then by Allard's Theorem, $\mathrm{spt}\,\nu_{-1}$ is the union of smooth disks near $x$, provided that $|\vec{A}|^2$ does not concentrate at $x$.

In §3, we improve the estimate (8) (hence (11)) slightly for the rescaled flows.

Then Theorem 1 follows in §4, because (11) bounds the number of concentration points of $|\vec{A}|^2$ in a given ball. Theorem 2 easily follows by a standard trick using an annulus. Such considerations of concentrated compactness originated for minimal surfaces in Schoen-Yau [SY] and Sacks-Uhlenbeck [SU]. We also mention a few examples and counterexamples relevant to Theorems 1 and 2.

Finally, in §5 we establish (slightly generalizing Ecker's Theorem) that if the quantity on the left-hand side of Theorem 4 is sufficiently small at all scales, then we have local regularity.

**Open Problems.**

**1.** Can the local Gauss-Bonnet estimate be proven via forms? Are there useful, analogous local estimates of various other characteristic classes?

**2.** *Conjecture:* If $M_0$ is embedded in $\mathbf{R}^3$, the rescaled surfaces converge smoothly (and with multiplicity one) to the blowup $N_t$. (That is, there are no concentration points or multiple layers.) Does the smooth convergence, at least, extend to $M_0$ immersed in $\mathbf{R}^3$ (see §4, Example 1)?

**3.** *General Conjecture:* The blowups are smooth for embedded hypersurfaces in $\mathbf{R}^n$ up to $n \leq 7$. What about immersed hypersurfaces? (Note: Velazquez [V] constructs a smooth, embedded surface evolving in $\mathbf{R}^8$ whose blowup is a minimizing cone, ergo not smooth.)

## §1 Local Gauss-Bonnet Estimate

Recall (see Fenchel [Fe]) that an immersed curve $\gamma$ in $\mathbf{R}^n$ satisfies

$$\text{(12)} \qquad \int_\gamma |\vec{k}| \geq 2\pi h(\gamma)$$

where $h(\gamma)$ is the number of components and $\vec{k}$ is the curvature of $\gamma$ in $\mathbf{R}^n$.

Suppose that $\gamma$ is the transverse intersection of manifolds $M^2$ and $N^{n-1}$ in $\mathbf{R}^n$. We want to estimate $|\vec{k}|^2$ in terms of $M$ and $N$. Let $\tau$, $n$ be an orthonormal frame for $T_x M$ along $\gamma$ such that $\tau$ is a unit tangent vector to $\gamma$. Let $e_1$ be a unit normal to $N$. Let $\vec{A}_M$, $\vec{A}_M$ be the second fundamental forms of $M$, $N$. The angle $\alpha$ between $T_x M$ and $T_x N$ satisfies

$$\sin^2 \alpha = (n \cdot e_1)^2.$$

We adopt the convention that if $P$ is a plane, then $P$ also denotes the orthogonal projection $\mathbf{R}^n \to P$.

**Lemma 5.**

$$\sin^2 \alpha |P \cdot \vec{k}|^2 + |P^\perp \cdot \vec{k}|^2 = |\vec{A}_M(\tau, \tau) - \vec{A}_N(\tau, \tau)|^2$$

where $P$ is any 2-dimensional vector space containing $n$ and $e_1$. (Note that $P$ is uniquely determined when $\sin^2 \alpha < 1$.)

**Proof.** Choose $e_2 \in T_x N \cap P$ so that $e_1, e_2$ is an orthonormal basis for $P$. Write
$$\vec{k} = (a, b, v), \qquad n = (n_1, n_2, 0),$$
according to the orthogonal decomposition $e_1$, $e_2$, $P^\perp$. Since $\vec{k} \cdot \tau = 0$, we have
$$\vec{A}_M(\tau, \tau) = T_x M \cdot \vec{k} = (n \cdot \vec{k})n = (n_1 a + n_2 b)(n_1, n_2, 0),$$
$$\vec{A}_N(\tau, \tau) = T_x N \cdot \vec{k} = (0, b, v).$$

so using $n_1^2 + n_2^2 = 1$,
$$|\vec{A}_M(\tau, \tau) - \vec{A}_N(\tau, \tau)|^2 = |(n_1^2 a + n_1 n_2 b, n_1 n_2 a + n_2^2 b - b, -v)|^2$$
$$= n_1^2 (n_1 a + n_2 b)^2 + n_1^2 (n_2 a - n_1 b)^2 + |v|^2$$
$$= n_1^2 (a^2 + b^2) + |v|^2,$$

as desired.

Putting $N = \partial B_r$, the Lemma yields

(13) $$|\vec{k}| \leq \frac{1}{\sin \alpha} \left( |\vec{A}_M| + \frac{1}{r} \right).$$

**Proof of Theorem 3.** Let $\phi = \phi(|x|) \in C_c^2(\mathbf{R}^n)$ with $\phi = 1$ on $B_1$, $\phi = 0$ on $\mathbf{R}^n \setminus B_R$, $\phi \geq 0$, $\phi' \leq 0$. Then by (9) and Fubini's Theorem,

$$\int_M \phi |\vec{A}|^2 d\mathcal{H}^2 = \int_0^1 \int_{M \cap \{\phi > u\}} |\vec{A}|^2 d\mathcal{H}^2 \, du$$
$$= \int_0^1 \left[ \int_{M \cap \{\phi > u\}} |\vec{H}|^2 d\mathcal{H}^2 - 8\pi c(u) + 8\pi g(u) + 4\pi h(u) \right.$$
$$\left. + \int_{M \cap \{\phi = u\}} 2\tilde{k} \cdot n \, ds \right] du$$

where $c(u)$ is the number of components of $M \cap \{\phi > u\}$, etc. Applying (12) together with the fact that $|\tilde{\vec{k}}| \leq |\vec{k}|$, we obtain

$$\begin{aligned}
\text{(last two terms)} &\equiv \int_0^1 4\pi h(u)\, du + \int_0^1 \int_{M \cap \{\phi = u\}} 2\tilde{\vec{k}} \cdot n \, ds \, du \\
&\leq \int_0^1 \int_{M \cap \{\phi = u\}} 4|\vec{k}|\, ds\, du \\
&\leq \int_0^1 \int_{M \cap \{\phi = u\}} \frac{4}{|\sin \alpha|} \left( \frac{1}{|x|} + |\vec{A}| \right) ds\, du \qquad \text{by (13)} \\
&= \int_M \frac{4|D_M \phi|}{|\sin \alpha|} \left( \frac{1}{|x|} + |\vec{A}| \right) d\mathcal{H}^2 \qquad \text{by co-area}
\end{aligned}$$

where $D_M \phi$ is $D\phi$ projected onto $T_x M$. Since $|D_M \phi| = |\sin \alpha| |D\phi|$, we obtain

$$\text{(last two terms)} \leq \int_M 4|D\phi| \left( \frac{1}{|x|} + |\vec{A}| \right) \leq \int_M \frac{4|D\phi|}{|x|} + \varepsilon \phi |\vec{A}|^2 + \frac{4|D\phi|^2}{\varepsilon \phi},$$

so

$$(14) \quad (1-\varepsilon) \int_M \phi |\vec{A}|^2 \leq \int_M \phi |\vec{H}|^2 + 8\pi g(M \cap B_R)$$
$$- 8\pi \int_0^1 c(u)\, du + \int_M \frac{4|D\phi|}{|x|} + \frac{4|D\phi|^2}{\varepsilon \phi}.$$

Set

$$\phi = \begin{cases} 1 & |x| \leq 1 \\ \dfrac{(R - |x|)^2}{(R-1)^2} & 1 \leq |x| \leq R \\ 0 & R \leq |x|, \end{cases}$$

which was selected so that $|D\phi|^2 / \phi = 4/(R-1)^2$. Calculate

$$\begin{aligned}
\int_M \frac{4|D\phi|}{|x|} &= \frac{8}{(R-1)^2} \int_{M \cap B_R \setminus B_1} \frac{R - |x|}{|x|} \\
&\leq \frac{8}{(R-1)^2} \int_0^{R-1} |\{x \in M : (R-|x|)/|x| > u\}|\, du \\
&\leq \frac{8}{(R-1)^2} \int_0^{R-1} \pi D' \left( \frac{R}{u+1} \right)^2 du \\
&= \frac{8\pi D' R}{R-1}.
\end{aligned}$$

Combining with (14) and using the fact that $c(M \cap B_S) \geq c'(M \cap B_R)$ for $S \in [1, R]$, we obtain Theorem 3.

**Proof of Theorem 4.** Take $R = 2$ in (10), scale to size $r$, observe that $g(M \cap B_R) \leq g(M_0)$, integrate with respect to time and use (7) and (8).

## §2 Blowups, Layers

We collect together a few facts concerning Brakke flows, blowups, layer structure, and local regularity.

A Radon measure $\mu$ is called *locally $k$-rectifiable* if the $k$-dimensional approximate tangent plane $T_x\mu$ exists $\mu$-a.e., that is to say, for $\mu$-a.e. $x$ the Radon measure $T_x\mu$ defined by

$$T_x\mu(A) = \lim_{\lambda \to 0} \lambda^{-k} \mu(x + A)$$

exists and is a positive multiple of $\mathcal{H}^k \lfloor P$ for some $k$-plane $P$. We will also use $T_x\mu$ to denote $P$ when no confusion arises. Equivalently, $\mu$ is the area measure $\mu = \mu_V = \|V\|$ of a rectifiable $k$-varifold $V = V_\mu$. The first variation of $\mu$ is defined by

$$\delta V_\mu(X) = \int \mathrm{div}_{S(x)} X(x) \, d\mu(x)$$

for $X \in C_c^\infty(\mathbf{R}^n, \mathbf{R}^n)$. By convention $S = S_\mu : \mathbf{R}^n \to G(n, k)$ denotes the $\mu$-measurable function that maps $x$ to the geometric tangent plane, denoted by $P$ above. Then $\mathrm{div}_S X = \sum_{i=1}^k D_{e_i} X \cdot e_i$, where $e_1, \ldots, e_k$ is any orthonormal basis of $S$. Recalling that $S$ also denotes the projection onto $S$, $\mathrm{div}_S X$ can be be written $S : DX$.

If the total first variation $\|\delta\mu_V\|$ is a Radon measure and is absolutely continuous with respect to $\mu$, then we define the generalized mean curvature vector $\vec{H} = \vec{H}_\mu \in L^1_{loc}(\mu)$ by requiring

(15) $$\int \mathrm{div}_S X \, d\mu = \int -\vec{H} \cdot X \, d\mu$$

for all $X \in C_c^\infty(\mathbf{R}^n, \mathbf{R}^n)$. Together with the previous equation, this is called the first variation formula. For more information about varifolds, see [S1].

A *Brakke flow* is a (measurable) family of Radon measures $\{\mu_t\}_{t \geq 0}$ on $\mathbf{R}^n$ satisfying equation (1) in the weak form

$$(16) \quad \int \phi \, d\mu_{t_2} \leq \int \phi \, d\mu_{t_1} + \int_{t_1}^{t_2} \int -\phi |\vec{H}|^2 + \vec{H} \cdot S^\perp \cdot D\phi \, d\mu_t \, dt.$$

for all nonnegative $\phi = \phi(x) \in C_c^1(\mathbf{R}^n)$ and $0 \leq t_1 \leq t_2$. If $\mu_t$ is the area measure of a smoothly moving surface $M_t$, then (16) is equivalent to (1), and holds with equality.

For nonsmooth flows, the inner integral is interpreted according to the following convention. If $\mu_t \llcorner \{\phi > 0\}$ is not a locally $k$-rectifiable Radon measure, the total variation $|\delta V_{\mu_t}| \llcorner \{\phi > 0\}$ is not absolutely continuous with respect to $\mu_t \llcorner \{\phi > 0\}$, or $\int \phi |\vec{H}|^2 \, d\mu = \infty$, then take the integral to be $-\infty$. Otherwise, all the quantities makes sense and the integral is taken literally. This definition is carefully arranged so that a Brakke flow $\mu_t$ is locally $k$-rectifiable for a.e. $t$ and the integral quantity is upper-semicontinuous with respect to weak convergence of Radon measures.

The formulation (16) is equivalent to [I3, 6.3], and slightly stronger than the original formulation of Brakke [B, 3.3]. For further details see Ilmanen [I3, 6.3].

The following lemma is of general interest. See also White [W].

**Lemma 7 (Monotonicity Formula for Brakke Flows).** For any Brakke flow $\{\mu_t\}_{t \geq 0}$ satisfying (6), any $(y, s) \in \mathbf{R}^n \times [0, \infty)$, and all $0 \leq t_1 \leq t_2 < s$,

$$\int \rho(x, t_2) d\mu_{t_2}(x) + \int_{t_1}^{t_2} \int \rho(x, t) \left| \vec{H}(x, t) + \frac{S^\perp(x, t) \cdot (x - y)}{2(s - t)} \right|^2 d\mu_t(x) \, dt$$
$$\leq \int \rho(x, t) \, d\mu_{t_1}(x)$$

where $\rho = \rho_{y,s}$, and $S(x,t) \equiv T_x \mu_t$. As in (16), the inner integral on the right-hand side is interpreted as $-\infty$ if one of the technical conditions fails.

**Proof.** 1. We adapt Huisken's proof to varifolds. By [B, 3.5], (16) implies that for any nonnegative test function $\phi = \phi(x,t) \in C_c^1(\mathbf{R}^n \times [0, \infty))$, $t_2 \geq t_1 \geq 0$,

$$(17) \quad \int \phi(\cdot, t_2) \, d\mu_{t_2}$$
$$\leq \int \phi(\cdot, t_1) \, d\mu_{t_1} + \int_{t_1}^{t_2} \int -\phi |\vec{H}|^2 + D\phi \cdot S \cdot \vec{H} + \frac{\partial}{\partial t} \phi \, d\mu_t \, dt.$$

Brakke [B, 5.8] proved the deep fact that that $S^\perp \cdot \vec{H} = \vec{H}$ $\mu$-a.e. whenever $|\delta V_\mu|$ is a Radon measure. Therefore whenever the inner integral is finite, we can calculate using the first variation formula (15)

$$\int -\phi |\vec{H}|^2 + D\phi \cdot S \cdot \vec{H} + \frac{\partial}{\partial t} \phi \, d\mu_t$$
$$\leq \int -\phi |\vec{H}|^2 + 2D\phi \cdot \vec{H} + S : D^2 \phi + \frac{\partial}{\partial t} \phi \, d\mu_t$$
$$= \int -\phi \left| \vec{H} - \frac{D\phi \cdot S^\perp}{\phi} \right|^2 + Q_S(\phi) \, d\mu_t$$

where
$$Q_S(\phi) = \frac{(D\phi \cdot S^\perp)^2}{\phi} + S : D^2 \phi + \frac{\partial}{\partial t} \phi.$$

When $\phi = \rho$, we obtain (remarkably)

$$Q_S(\phi) = 0$$

for any $k$-plane $S$.

2. To insert $\rho$ into the above formula, let $\psi = \psi_R$ be a cutoff function with $\chi_{B_R(y)} \leq \psi \leq \chi_{B_{2R}(y)}$, $R|D\psi| + R^2 |D^2 \psi| \leq C$. Calculate

$$Q_S(\psi \rho) = \psi Q_S(\rho) + \rho Q_S(\psi) + 2 D\psi \cdot D\rho \leq 0 + C \left( \frac{1}{R^2} + \frac{1}{s-t} \right) \chi_{B_{2R}(y) \setminus B_R(y)} \rho$$

where we used the fact that $|D\rho| \leq |x-y|\rho/2(s-t)$. Inserting $\psi\rho$ above, we obtain

$$\int \psi\rho \, d\mu_{t_2} + \int \psi\rho \left| \vec{H} + \frac{S^\perp \cdot (x-y)}{2(s-t)} + \frac{D\psi}{\psi} \right|^2$$
$$\leq \int \psi\rho \, d\mu_{t_1} + \left( \frac{1}{R^2} + \frac{C}{s-t_2} \right) \int_{t_1}^{t_2} \int_{B_{2R}(y) \setminus B_R(y)} \rho \, d\mu_t \, dt.$$

By the monotone convergence theorem and the dominated convergence theorem, this will establish the desired result, once we prove that $\sup_{[t_1,t_2]} \int \rho < \infty$. This is a quick calculation: insert the function $\phi = \eta_R(|x-y|^2 + 2kt)$ into (16), where

$$\eta_R(u) = \begin{cases} (R^2 - u)^4 & u \leq R^2 \\ 0 & u > R^2. \end{cases}$$

As shown in [B, 3.7], $\int \phi \, d\mu_t$ is nonincreasing. Using (6), it follows that

$$\mu_t(B_R(y)) \leq 4 \cdot 2^k \omega_k D R^k, \qquad t \in [0, R^2/4k],$$

for all $y \in \mathbf{R}^n$, $R > 0$. In particular this holds for $t \in [0, t_2]$ provided that $R^2 \geq t_2/C$. From this we may calculate that

$$\int \rho_{y,s} \, d\mu_t \leq CD,$$

where $C$ is bounded for fixed $s$ and $t \in [t_1, t_2]$. This completes the proof of Lemma 7.

Next, define the parabolically recaled measures at a point $(y,s)$ by

$$\mu_t^{(y,s),\lambda}(A) = \lambda^{-k} \mu_{s+\lambda^2 t}(y + \lambda \cdot A), \qquad t \in [-s/\lambda^2, \infty).$$

For the following lemma see also White [W].

**Lemma 8 (Weak Existence of Blowups).** For any Brakke flow satisfying (6), any point $(y,s)$, and any sequence $\lambda_i \to 0$ there exists a subsequence $\{\lambda_{i_j}\}_{j\geq 1}$ and a limiting Brakke flow $\{\nu_t\}_{t\in\mathbf{R}}$ such that $\mu_t^{(y,s),\lambda_{i_j}} \rightharpoonup \nu_t$ in the sense of Radon measures for all $t \in \mathbf{R}$ and

(a) $$\nu_t(A) = \nu_t^\lambda(A) \equiv \lambda^{-k}\nu_{\lambda^2 t}(\lambda \cdot A), \qquad t < 0,$$

for all $\lambda > 0$, and $\nu_{-1}$ satisfies

(b) $$\vec{H}(x) + \frac{S(x)^\perp \cdot x}{2} = 0, \qquad \nu_{-1}\text{-a.e. } x.$$

Furthermore, Huisken's integral converges:

(c) $$\int \rho_{0,0}(x,-1)\,d\nu_t(x) = \lim_{t'\uparrow s}\int \rho_{y,s}\,d\mu_{t'}, \qquad t < 0.$$

**Proof.** 1. Write $\mu_t^\lambda = \mu_t^{(y,s),\lambda}$. By (6) and Lemma 7, we have

(18) $$\mu_t^\lambda(B_R(x)) \leq C(k)DR^k, \qquad x \in \mathbf{R}^k,\quad t \geq 0,\quad R > 0.$$

Then by the Compactness Theorem [I3, 7.1], there exists a subsequence $\lambda_{i_j} \to 0$ and a limit Brakke flow $\{\nu_t\}_{t\in\mathbf{R}}$ such that $\mu_t^{\lambda_{i_j}} \rightharpoonup \nu_t$ for $t \in \mathbf{R}$. Then

$$\int \rho_{0,0}\,d\mu_t \leq \liminf_{j\to 0}\int \rho_{0,0}\,d\mu_t^{\lambda_{i_j}}.$$

Using (18) and the exponential decay of $\rho$ we calculate for $t < 0$,

$$\int_{\mathbf{R}^n\setminus B_R}\rho_{0,0}\,d\mu_t^{\lambda_i} \leq \frac{C(k)}{(-t)^{k/2}}\sum_{j\geq 1}\int_{B_{R^{j+1}}\setminus B_{R^j}}e^{-R^{2j}/(-4t)}\,d\mu_t^{\lambda_i}$$
$$\leq \frac{C(k)D}{(-t)^{k/2}}\sum_{j\geq 1}R^{(j+1)k}e^{-R^{2j}/(-4t)}$$
$$= Df(R^2/(-t))$$

where $f(R^2/(-t))$ is a function that converges to 0 as $R \to \infty$, for each fixed $t$. Therefore the sequence of measures $\rho_{0,0}\,d\mu_t^{\lambda_i}$ is tight for each $t < 0$ and

$$\int \rho_{0,0}\,d\nu_t = \lim_{j\to 0}\int \rho_{0,0}\,d\mu_t^{\lambda_{i_j}} = \lim_{t'\uparrow s}\int \rho_{y,s}\,d\mu_t', \qquad t < 0,$$

which is (c).

2. Then by Lemma 7 applied to $\{\nu_t\}_{t \in \mathbf{R}}$, for a.e. $t < 0$ $\nu_t$ is locally $k$-rectifiable with $|\vec{H}| \in L^2_{loc}(\mu)$, and

$$\vec{H} + \frac{S^\perp \cdot x}{-2t} = 0, \qquad \nu_t\text{-a.e. } x. \tag{19}$$

Next we show self-similarity. Define $\tilde{\nu}_t(A) = (-t)^{-k/2}\nu_t((-t)^{1/2}A)$, $t < 0$. It suffieces to show that $\tilde{\nu}_t$ is constant. By (17), calculate for $t_1 \leq t_2 < 0$, $\phi \in C^2_c(\mathbf{R}^n, [0, \infty))$,

$$\int \phi \, d\tilde{\nu}_{t_2} - \int \phi \, d\tilde{\nu}_{t_1}$$

$$= \frac{1}{(-t_2)^{k/2}} \int \phi((-t_2)^{1/2}x) \, d\nu_{t_2} - \frac{1}{(-t_1)^{k/2}} \int \phi((-t_1)^{1/2}x) \, d\nu_{t_1}$$

$$\leq \int_{t_1}^{t_2} \frac{k}{2(-t)^{k/2+1}} \int \phi \, d\nu_t + \frac{1}{(-t)^{k/2}} \int -\phi|\vec{H}|^2 + (-t)^{1/2} D\phi \cdot S^\perp \cdot \vec{H}$$

$$- \frac{1}{2(-t)^{1/2}} D\phi \cdot x \, d\nu_t \, dt$$

$$= \int_{t_1}^{t_2} \frac{1}{(-t)^{k/2}} \int -\frac{k}{2t}\phi - \frac{\phi \vec{H} \cdot S^\perp \cdot x}{2t} + (-t)^{1/2} D\phi \cdot S^\perp \cdot \frac{x}{2t}$$

$$- \frac{1}{2(-t)^{1/2}} D\phi \cdot x \, d\nu_t \, dt$$

after substituting from (b). From the first variation formula and Brakke [B, 5.8] we have

$$\int -\frac{\phi \vec{H} \cdot S^\perp \cdot x}{2t} \, d\nu_t = \int -\frac{\phi \vec{H} \cdot x}{2t} \, d\nu_t = \int \frac{\mathrm{div}_S(\phi x)}{2t} \, d\nu_t$$

$$= \int \frac{k\phi}{2t} + (-t)^{1/2} \frac{D\phi \cdot S \cdot x}{2t} \, d\nu_t.$$

Combining with the above, this implies that $\int \phi \, d\tilde{\nu}_t$ is nonincreasing in $t$.

Next, assume without loss of generality that $\phi < \rho_{0,0}$ and apply the same calculation to $\rho_{0,0} - \phi$ (using the exponential decay of $\rho_{0,0}$ and (18) to validate the insertion of this function) to show that $\int \rho_{0,0} - \phi \, d\tilde{\nu}_t$ is also nonincreasing.

It follows by (c) that $\int \phi \, d\tilde{\nu}_t$ is constant, that is, (a) holds. Then by (19), (b) holds. This completes the proof of Lemma 8.

The next Theorem follows from [S1, 24.2].

**Theorem 9 (Allard's Regularity Theorem).** There exist $\varepsilon_1 = \varepsilon_1(n,k)$, $\sigma = \sigma(n,k) > 0$ with the following property. Let $\mu$ be an integer $k$-rectifiable Radon measure such that $|\vec{H}| \in L^1_{loc}(\mu)$. If $\varepsilon \in (0, \varepsilon_1]$, $0 \in \operatorname{spt} \mu$, $K > 0$, $r > 0$,

$$|\vec{H}| \leq \frac{\varepsilon}{r}, \qquad \mu\text{-a.e. } x \in B_r,$$

and

$$\mu(B_r) \leq (1+\varepsilon)\omega_k r^k,$$

then there is a $k$-plane $T \ni 0$, and domain $\Omega \subseteq T$, and a $C^{1,\alpha}$ vector function $u : \Omega \to T^\perp$ such that

$$\operatorname{spt} \mu \cap B_{\sigma r} = \operatorname{graph} u \cap B_{\sigma r}$$

where $\operatorname{graph}(u) = \{x + u(x) : x \in \Omega\}$, and

$$\sup \left|\frac{u}{r}\right| + \sup |Du| + r^\alpha [Du]_\alpha \leq C\varepsilon^{1/4n}$$

where $[Du]_\alpha = \sup |Du(x) - Du(y)|/|x - y|^\alpha$.

To obtain the principal hypothesis above (which we call the single-layer hypothesis) we use the following.

**Theorem 10 (Simon's Lemma on $|\vec{A}|^2$ [S2, Lemma 4]).** For each $n \geq 3$, $D > 0$, there is $\varepsilon_2 = \varepsilon_2(n, D)$ such that if $M$ is a smooth 2-manifold properly immersed in $B_R$ and

$$\int_{M \cap B_R} |\vec{A}|^2 \leq \varepsilon^2 \leq \varepsilon_2^2, \qquad \mathcal{H}^2(M \cap B_R) \leq D\pi R^2$$

then there are pairwise disjoint closed disks $\bar{P}_1, \ldots, \bar{P}_N$ in $M \cap B_R$ ($B_R$ open) such that

$$\sum_m \operatorname{diam}(P_m) \leq C(n, D)\varepsilon^{1/2} R$$

and for any $S \in [R/4, R/2]$ such that $M$ is transverse to $B_S$ and $\partial B_S \cap \cup_m P_m = \emptyset$, we have

$$M \cap B_S = \cup_{l=1}^m D_l$$

where each $D_i$ is an *embedded* disk. Furthermore, for each $D_l$ there is a 2-plane $L_l \subseteq \mathbf{R}^n$, a simply connected domain $\Omega_l \subseteq L_l$, disjoint closed balls $\bar{B}_{l,p} \subseteq \Omega_l$, $p = 1, \ldots, p_l$ and a function

$$u_l : \Omega_l \setminus \cup \bar{B}_{l,p} \to L_l^\perp$$

such that

$$\sup \left| \frac{u_l}{R} \right| + |Du_l| \leq C(n, D)\varepsilon^\gamma$$

where $\gamma = 1/2(2n - 3)$, and

$$D_l \setminus \cup_m \bar{P}_m = \operatorname{graph}(u_l|\Omega_l \setminus \cup_p \bar{B}_{l,p}).$$

The above Theorem is stated in [S2] for embedded surfaces but is valid for immersed ones. By applying it to subballs of $B_R$ at various scales, we obtain

**Corollary 11.** *Under the above hypotheses, for any $x \in M \cap B_R$ such that $B_{2r}(x) \subseteq B_R$, the connected component $M'$ of $M \cap B_r(x)$ containing $x$ is embedded and satisfies*

$$\mathcal{H}^2(M') \leq \mathcal{H}^2(\operatorname{graph}(u_l|\Omega')) + \mathcal{H}^2(\cup_m P_m) \leq \pi r^2(1 + C(n, D)\varepsilon^\gamma),$$

$$\mathcal{H}^2(M') \geq \pi r^2(1 - C(n, D)\varepsilon^{2\gamma}),$$

where $\Omega'$ is an appropriate domain in a 2-plane.

## §3 Estimate of H on a Blowup Sequence

We will establish an improved estimate for $\iint |\vec{H}|^2$ on a blowup sequence. Integrating (2) and using (6), we have

$$(20) \qquad \int_{-2}^{0} \int_{M_t^\lambda} \rho_{0,0}(x,t) \left| \vec{H} + \frac{S^\perp \cdot x}{-2t} \right|^2 d\mathcal{H}^k(x)\, dt \equiv \delta(\lambda) \to 0$$

as $\lambda \to 0$. Using the triangle inequality and (7), we have for any ball $B_r(x) \subseteq B_R$,

$$\int_{-1-\tau}^{-1} \int_{M_t^\lambda \cap B_r(x)} |\vec{H}|^2 \leq C(k) D\tau r^k R^2 + \delta_R(\lambda),$$

where $\delta_R(\lambda) \to 0$ as $\lambda \to 0$. This estimate expresses that $|\vec{H}| \leq |x| \in L^\infty_{loc}$ in the limit, coming from (5).

We improve on this slightly as follows.

**Lemma 6.** Let $\{M_t\}_{t \in [0,T)}$ satisfy (1) and (6) and let $\{M_t^\lambda\}_{t \in [-T/\lambda^2, 0)}$ be a blowup sequence. Then if $B_{2r}(x) \subseteq B_R$, we have

$$\int_{-1-\tau}^{-1} \int_{M_t^\lambda \cap B_r(x)} |\vec{H}|^2 \leq C(k) D\tau (r^k + r^{k-1} R) + \delta_R(\lambda)$$

where $\delta_R(\lambda) \to 0$ as $\lambda \to 0$.

**Note.** For self-shrinkers in $\mathbf{R}^3$, this leads to $\int_{M_{-1} \cap B_R} |\vec{H}|^2 \leq CDR^2$. Compare to known examples:

$$\begin{aligned} \text{Compact:} &\quad C \\ \text{Cylinder:} &\quad CR \\ \text{Conelike at } \infty: &\quad C \log R. \end{aligned}$$

For examples of the latter see [C], [ACI]. We believe that all examples fall into these classes and therefore Lemma 6 is not optimal.

**Proof.** Let $\zeta$ be a smooth function of compact support, $M$ an immersed surface. Using the first variation formula (15),

$$\int_M \zeta^2 |\vec{H}|^2 = \int_M \zeta^2 \vec{H} \cdot \left(\vec{H} + \frac{S^\perp \cdot x}{-2t}\right) + \zeta^2 \frac{\vec{H} \cdot x}{2t}$$

$$\leq \int_M \frac{1}{2}\zeta^2|\vec{H}|^2 + \frac{1}{2}\zeta^2 \left|\vec{H} + \frac{S^\perp \cdot x}{-2t}\right|^2 - \frac{1}{2t}\text{div}_M(\zeta^2 x)$$

so taking $\zeta$ to be a standard cutoff function for $B_r(x)$ inside $B_{2r}(x)$ and using (7),

$$\int_{M \cap B_r(x)} |\vec{H}|^2 \leq \int_{M \cap B_{2r}(x)} \left|\vec{H} + \frac{S^\perp \cdot x}{-2t}\right|^2 + \frac{C(k)Dr^k}{-2t}\left(1 + \frac{R}{r}\right).$$

Let $M = M_t^\lambda$, integrate with respect to $t$, and recall (20).

# §4 Partial Regularity

Let $\{M_t\}_{t \in [0,T)}$ be a family of 2-dimensional surfaces of finite genus properly immersed in $\mathbf{R}^n$ that satisfy (1), (3), (6). Suppose $\{\mu_t^{\lambda_j}\}_{t \in [-T/\lambda_j^2, 0)}$ is a blowup sequence converging to $\{\nu_t\}_{t<0}$ as in Lemma 8. Write $N_t = \text{spt } \nu_t$, $N = N_{-1}$, $\nu = \nu_{-1}$.

**Theorem 1 (more precisely).** In the above situation, there is a discrete set $Q \subseteq \mathbf{R}^n$ such that

(i) $\qquad \#(Q \cap B_R) \leq C(n, D)[R^2 + 1 + 16\pi g(M_0)], \qquad R > 0,$

and a smooth 2-manifold $X$ equipped with a proper, surjective, conformal branched immersion $\iota : X \to N$, branched over $Q$, solving (5) in the parametrized sense, such that

(ii) $\qquad \#\iota^{-1}(x) = \theta^2(\nu, x), \qquad x \in N \setminus Q,$

where $\theta^2(\nu, x) = \lim_{r \to 0} \nu(B_r(x))$, and sequences $j_k \to \infty$, $t_k \to 1$, $R_k \to \infty$, and $r_{k,q} \to 0$ for each $q \in Q$, and a diffeomorphism

(iii) $\quad X \cap \iota^{-1}(B_{R_k} \setminus Q) \cong Y_k \equiv \eta_k^{-1}(B_{R_k} \setminus \cup_{q \in Q} \bar{B}_{r_{q,k}}(q)) \subseteq M_{t_{j_k}}^{\lambda_{j_k}}$

where $\eta_k$ is the map immersing $M_{t_{j_k}}^{\lambda_{j_k}}$ in $\mathbf{R}^n$, and $Y_k$ has a smooth boundary.

In consequence $g(X) \leq g(M_0)$, though there is presently no bound on the number of ends of $X$.

For the definition of conformal branched immersion solving (5), see Gulliver [G1, p.276]. In particular, $\iota$ restricted to $X \setminus \iota^{-1}(Q)$ is an immersion.

A *branch point* is a point $y \in X$ such that $D\iota(y)$ is not an injection. The set of branch points is discrete. If $\iota(B_r(y))$ is smooth for some $r > 0$, we call $y$ a *false branch point*. Otherwise $y$ is a *true branch point*.

The image of the branch points lies in $Q$. The image of the false branch points may depend on the subsequence $\{j_k\}$. However, by the structure near a true branch point [G1], the image of the true branch points lies in a subset $P$ of $Q$ that can be characterized by the geometry of $N$ alone, and is therefore independent of the subsequence $\{j_k\}$.

In a later version of this paper, we will show that $Q$ is finite and that $\int |\vec{A}|^2$ concentrates only in multiples of $4\pi$.

**Proof of Theorem 1.** 1. Define $\sigma_j = |\vec{A}|^2 \mathcal{H}^2 \lfloor M_{t_j}^{\lambda_j}$. First we examine the concentration behavior of $\sigma_j$ as $j \to \infty$. From Theorem 3 (with $\varepsilon = 1/2$, $R = 2$) and Lemma 6 and recalling (7), we deduce

$$\frac{1}{2} \int_{-1-\tau}^{-1} \int_{M_t^{\lambda_j} \cap B_R} |\vec{A}|^2 \leq \int_{-1-\tau}^{-1} \int_{M_t^{\lambda_j} \cap B_{2R}} |\vec{H}|^2 + \tau(8\pi g(M_0) + CD)$$
$$\leq \tau(CDR^2 + 8\pi g(M_0) + CD) + \delta_R(j),$$

and

$$\int_{-1-\tau}^{-1} \int_{M_t^{\lambda_j}} \rho_{0,0} \left| \vec{H} + \frac{S^\perp \cdot x}{-2t} \right|^2 \leq \delta(j).$$

Select $\tau_j \to 0$, $R_j \to \infty$ so that $\delta_R(j)/\tau_j$, $\delta(j)/\tau_j \to 0$. Then choose $t_j \in [-1 - \tau_j, -1]$ such that

$$\sigma_j(B_R) \leq CD(R^2 + 1) + 16\pi g(M_0) + 2\delta_R(j)/\tau_j,$$

and

(21)
$$\int_{M_{t_j}^{\lambda_j}} \rho_{0,0} \left| \vec{H} + \frac{x \cdot S^\perp}{-2t_j} \right|^2 \leq \delta(j)/\tau_j.$$

Select a subsequence $\{j_k\}$ such that $\sigma_{j_k}$ converges to a Radon measure $\sigma$. Using Lemma 8, we see that $\mu_{t_{j_k}}^{j_k} \rightharpoonup \nu$ as Radon measures. Applying the monotonicity formula at points other than $(0,0)$, together with Lemma 8(a), it follows that the supports converge locally in the sense of Hausdorff distance.

2. Define $\varepsilon_3 = \varepsilon_3(n, D) = \min(\varepsilon_2, \varepsilon_1, (\varepsilon_1/C(n,D))^{1/\gamma})$, where $\varepsilon_1$, $\varepsilon_2$, $\gamma$ are as in Lemmas 9 and 10, and $C(n, D)$ will emerge below. Define the set $Q$ of *concentration points* to consist of $p \in N$ such that $\sigma(p) \geq \varepsilon_3$. Then $Q$ is discrete and

$$\#(Q \cap B_R) \leq \frac{CD(R^2 + 1) + 16\pi g(M_0)}{\varepsilon_3}$$

which gives (i).

Suppose $\varepsilon \in (0, \varepsilon_3]$, $R, r > 0$, $p \in N$ and

(22) $\qquad B_r(p) \subseteq B_R, \qquad r \in (0, \varepsilon/R], \qquad \sigma(B_r(p)) < \varepsilon^2.$

(Note that this can be arranged for any $p \in N \setminus Q$ by choosing $r$ sufficiently small.) We claim that $N$ is the image of an immersed manifold near $p$. Write $M_k = M_{t_{j_k}}^{\lambda_{j_k}}$, $\mu_k = \mathcal{H}^2 \lfloor M_k$. Note that $\sigma_{j_k}(B_r(p)) < \varepsilon^2$ for sufficiently large $k$. Then by Theorem 10 (Simon) and Corollary 11, there are $s_k \in [r/4, r/2]$, and a decomposition of $M_k \cap B_{s_k}(p)$ into distinct components $D_{k,l}$, each of which is an embedded disk satisfying

$$\mathcal{H}^2(D_{k,l} \cap B_\rho(x)) \leq (1 + C(n,D)\varepsilon^\gamma)\rho^2, \qquad x \in D_{k,l}, \quad |x-p| + \rho \leq s_k.$$

By choosing a further subsequence of $\{j_k\}$, labelled the same, we may assume that as $k \to \infty$,

$$l_k \to l_0, \qquad s_k \to s \in [r/4, r/2], \qquad \mathcal{H}^2 \lfloor D_{k,l} \to \nu_l, \qquad l = 1, \ldots, l_0.$$

Note by (7), (8), and Allard's Compactness Theorem [S1, 42.7] that $V_{\mathcal{H}^2 \lfloor D_{k,l}} \to V_{\nu_l}$ as varifolds. It follows by the lower semicontinuity of the quantity in (21) (easily checked, see [I3] for an analogous calculation) that $\nu_l$ weakly solves (5).

Then the hypothesis on $r$ implies that $|\vec{H}| \leq \varepsilon/r$. Also

$$\nu_l(B_\rho(x)) \leq (1 + C(n,k)\varepsilon^\gamma)\rho^2, \qquad x \in \mathrm{spt}\,\nu_l \cap B_s(p), \quad |x-p| + \rho \leq s.$$

By the choice of $\varepsilon_3$, Theorem 9 (Allard) applies. Also employing the Schauder estimates, $\mathrm{spt}\,\nu_l \cap B_{\sigma s}(p)$ is the graph of a $C^\infty$ function $u$ defined over a domain in a 2-plane $L_l$ with the estimate

$$(23) \qquad \left|\frac{u}{r}\right| + \sup|Du| + r\sup|Du|^2 \leq C(n,D)\varepsilon^{\gamma_1}$$

where $\gamma_1 = \gamma/4n = 1/8n(n-3)$. Also $\sum \nu_l = \nu \lfloor B_s$, so $N \cap B_{\sigma s}(p) = \cup_l \mathrm{spt}\,\nu_l$ is the image of an immersed manifold.

3. Next we construct a 2-manifold $Y$ and a surjective immersion $\iota : Y \to N \setminus Q$ satisfying (ii).

Select balls $\{B_{r(x)}(x)\}_{x \in N \setminus Q}$ satisfying (22). By taking a subcover and diagonalizing, there exist points $\{x_i\}_{i \geq 1} \subseteq N \setminus Q$ and radii $s_i$ such that $\{B_{\sigma s_i}(x_i)\}_{i \geq 1}$ is a locally finite covering of $N \setminus Q$, together with sequences $s_{k,i} \to s_i$, $D_{k,i,l} \to D_{i,l}$ (in Hausdorff distance) where $D_{i,l} \subseteq B_{\sigma s_i}(x_i)$ satisfies (23). Applying (22) and its consequences to subballs of $B_{r_i}$, we may assume that $D_{i,l} \subseteq D'_{i,l} \subseteq N$, where $D'_{i,l}$ is a smooth, embedded disk with boundary in $\partial B_{2\sigma s_i/2}(x_i)$.

Let $U$ be an open set containing $D'_{i,l}$ on which the nearest point projection $\pi_{i,l} : U \to D'_{i,l}$ is well-defined. For sufficiently large $k$ (depending on $i$), $D_{k,i,l} \subseteq U$.

If $D_{k,i,l}$ overlaps $D_{k,i',l'}$ (that is, in the domain of the immersion defining $M_k$) for infinitely many indices $k$, then $D'_{i,l}$ overlaps $D'_{i',l'}$ (in $\mathbf{R}^n$) and $D'_{i,l} \cup D'_{i',l'}$ is a smooth, embedded surface with estimates like (23). Then for $k$ sufficiently large (depending on $i$, $i'$),

$$(24) \qquad \pi_{i,l}(x) = \pi_{i',l'}(x), \qquad x \in D_{k,i,l} \cup D_{k,i',l'}.$$

Let $P_{k,i,m}$ denote the bad disks in $D_{k,i,l}$. By Lemma 10, we can select monotone sequences $R_k \to \infty$ and $r_{k,q} \to 0$ for each $q \in Q$, such that the set $F_k \equiv B_{R_k} \setminus \cup_{q \in Q} B_{r_{k,q}}(q)$ satisfies

(a) $\partial F_k \cap (\cup P_{k,i,m}) = \emptyset$ and $\partial F_k$ is transverse to $M_k$, so that $Y_k$ defined as $M_k \cap F_k$ is a compact, smooth manifold with boundary, and

(b) Condition (24) holds for all $B_{r_i}(x_i)$, $B_{r'_i}(x'_i)$ that meet $F_k$.

Then the map $\pi_k : Y_k \to N \setminus Q$ characterized by

$$\pi_k(x) = \pi_{D_{i,l}}(x), \qquad x \in D_{k,i,l},$$

is well-defined.

Now define a homotopic map $\iota_k : Y_k \to N \setminus Q$ that equals $\pi_k$ on $Y_k \setminus \cup P_{k,i,m}$, and such that $\iota_k | P_{k,i,m} : P_{k,i,m} \to \pi_k(P_{k,i,m})$ is a diffeomorphism. This is possible since each $P_{k,i,m}$ is a disk and $P_{k,i,m}$ is disjoint from $\partial F_k \cup (\cup_i \partial B_{\beta s_i}(x_i))$.

Then $\iota_k$ is an immersion into $N \setminus Q$. By an additional small homotopy, we may arrange that $\iota_k(Y_k) = N \cap F_k$, $\iota_k(\partial Y_k) = N \cap \partial F_k$. Note that for each $x \in N \setminus Q$, $\#\iota_k^{-1}(x) = \theta^2(\mu, x)$ for sufficiently large $k$.

We endow $Y_k$ with the metric $g_k$ induced by $\iota_k$. Choosing a further subseqence of $\{j_k\}$, also labelled $\{j_k\}$, we pass $(Y_k, g_k, \iota_k)$ to limits, obtaining a smooth 2-manifold $(Y, g)$ and a surjective immersion $\iota : Y \to N \setminus Q$ satisfying (ii).

4. In fact for each compact $K \subseteq \mathbf{R}^n \setminus Q$, the manifolds $\iota_k^{-1}(K) \subset Y_k$ are eventually isometric. Then by redefining $R_k$, $r_{k,q}$ to converge more slowly, we may arrange isometric embeddings

$$Y_1 \hookrightarrow Y_2 \hookrightarrow \cdots \qquad \hookrightarrow Y$$

so $Y_k \cong \iota^{-1}(F_k) \subseteq Y$ (isometric).

By Lemma 12 below, for each $q \in Q$ there is $r > 0$ such that

$$\iota^{-1}(B_r(q) \setminus \{q\}) = \cup_h A_h$$

where each $A_h$ is an annulus. Thus by making $R_k$ converge even more slowly, we may arrange that $\iota^{-1}(B_{r_k}(q) \setminus \{q\}) \subseteq Y$ consists of annuli for all $B_{r_k}(q)$ that meet $B_{R_k}$, so $\iota^{-1}(F_k) \cong \iota^{-1}(B_{R_k} \setminus Q)$ (diffeomeorphic). This yields conclusion (iii).

Define $X$ to be $Y \cup \{a_1, a_2, \ldots\}$ where $a_h$ is added to compactify $A_h$ along the open boundary component that maps to $q$. Then $\iota$ extends continuously to $X$ by defining $\iota(a_h) = q$. By Gulliver [G2], the conformal structure (though not the metric) of $Y$ extends to $X$ in such a way that $\iota$ satisfies (5) with respect to this conformal structure. Therefore $\iota$ is a conformal branched immersion. This completes the proof of the Theorem, subject to the Lemma.

**Lemma 12.** There is $\varepsilon_4 = \varepsilon_4(n, D)$ such that if $\varepsilon \in (0, \varepsilon_4]$, $R, r > 0$, and

$$B_r(q) \subseteq B_R, \qquad r \in (0, \varepsilon/R], \qquad \sigma(B_r(q) \setminus \{q\}) < \varepsilon^2,$$

then
$$\iota^{-1}(B_{r/2}(q) \setminus \{q\}) = (\cup_e D_e) \cup (\cup_h A_h)$$

where each $D_e$ is a smooth embedded disk and each $A_h$ is an annulus, properly immersed in $B_{r/2}(q) \setminus \{q\}$ with a well-defined branch index $\kappa$ about $q$, and
$$\left|\mathcal{H}^2(A_h) - \kappa\pi\left(\frac{r}{2}\right)^2\right| \leq C(n,D)\varepsilon^{\gamma_1/2}r^2.$$

**Proof.** We return to the suppression of $\iota$. Let $\beta \in (0, 1/20]$. Let $A$ be a connected component of $X \cap B_{r/2}(q) \setminus \{q\}$. Let $E$ be a connected component of $A \setminus B_{\beta r}(q)$. If $E$ is disjoint from $B_{\beta r}(q)$, then by (23), $E$ is a disk, so $A = E$ and we are done.

Fix $z \in E \cap \partial B_{r/4}$ and let $L = z + T_z E$. Note by (23) that $|\vec{A}| \leq C\varepsilon^{\gamma_1}/\beta r$ on $E$, so
$$|\nu(x) - \nu(z)| \leq \frac{C(n,D)K\varepsilon^{\gamma_1}}{\beta}, \qquad x \in E_K,$$

where $E_K = \{x \in E : d_E(x,z) < Kr\}$, where $d_E$ is the intrinsic distance in $E$. Assume $\varepsilon$ is small enough (depending on $K$, $\beta$) to ensure

(25) $$\frac{C(n,D)K\varepsilon^{\gamma_1}}{\beta} \leq \beta.$$

Then by (23) $L : E_K \to L$ is an immersion with a local inverse $u$ (near each point) satisfying
$$\left|\frac{u}{r}\right| + |Du| \leq \beta, \qquad x \in E_K.$$

Since $E$ meets $B_{\beta r}(q)$, it follows that $L$ meets $B_{2\beta r}$, and $E_L$ lies within $3\beta r$ of $T = q + T_z E$. Therefore $E' \equiv E_K \setminus B_{4\beta r}$ consists of a fairly flat ribbon that winds around $B_{4\beta r}$ at least $K/(2\pi+1)$ times. In particular, if $E'$ does not meet itself coming around, then $E'$ contains a $K/(2\pi+1)$-valued graph over $T \cap B_{(1/2-3\beta)r}(q) \setminus B_{4\beta r}(q)$ and
$$\frac{K}{2\pi+1}[(1/2 - 3\beta)^2 - (4\beta)^2]r^2 \leq \mathcal{H}^2(E') \leq CDr^2.$$

Select $K = (2\pi + 1)(200/33)CD + 1$ so that this is impossible. We conclude that $E'$ does meet itself, and increasing $K$ by 1, $E' = E \setminus B_{4\beta r}$. In particular $E \setminus B_{4\beta r}$ is an immersed annulus with one boundary component lying in $\partial B_{4\beta r}(q)$. Iterating this construction inward to $q$, we establish that $A$ is an annulus (or a disk disjoint from $q$), and $E = A \setminus B_{\beta r}(q)$. By Gulliver [G2], the annulus $A$ is conformally equivalent to the punctured disk $B_1^2 \setminus \{0\}$, and the immersion that parametrizes $A$ can be extended smoothly to a branched immersion branched over $q$. The value of $\kappa$ is the same at all scales, and equal to the branching index.

This holds subject to (25). To activate (25), select $K$ as above, and $\beta = (C(n,D)K\varepsilon^{\gamma_1})^{1/2}$ so as to satisfy (25). Choose $\varepsilon_4$ so that $\varepsilon \leq \varepsilon_4$ ensures that $\beta \leq 1/20$.

Finally,

$$\left| \mathcal{H}^2(A) - \kappa\pi \left(\frac{r}{2}\right)^2 \right| \leq |\mathcal{H}^2(E') - \mathcal{H}^2(T \cap B_{r/2} \setminus B_{4\beta r})| + \mathcal{H}^2(A \cap B_{4\beta r})$$
$$\leq C(n,D)\beta r^2 \leq C(n,D)\varepsilon^{\gamma_1/2} r^2.$$

This completes the proof of Lemma 12.

**Proof of Theorem 2.** Suppose $n = 3$ and $M_0$ is embedded.

1. Let $q \in N \setminus Q$. Since $M_0$ is embedded, each $M_t^{\lambda_j}$ is embedded, so the disks $D_{k,l}$ of the previous proof are disjoint for each fixed $k$. Passing to limits, we find that $\operatorname{spt} \nu_l$ is the smooth boundary in $B_s(q)$ of an open set $E_l \subseteq B_s(q)$ such that for every $l'$, either $E_l \subseteq E_{l'}$ or $E_{l'} \subseteq E_l$. Then by the strong maximum principle, each pair $\operatorname{spt} \nu_l$, $\operatorname{spt} \nu_{l'}$ is either coincident or disjoint. Therefore, $N \setminus Q$ is embedded.

2. Select $r$ small enough that Lemma 12 applies. Let $A$ be any annulus component of $X \cap B_{r/2}(q) \setminus \{q\}$. Then $\iota(A)$ is smooth, properly embedded in $B_{r/2}(q)$, and solves (5) away from $q$. Lemma 12 applies to the embedding $\iota(A) \to \mathbf{R}^n$ with $\kappa = 1$ to yield the area bound needed for Allard's Theorem. By (7) and an elementary cutoff function argument, the Radon measure $\mathcal{H}^2 \lfloor \iota(A)$ satisfies (5) in the weak sense in $B_{r/2}(q)$. Then by Allard's Theorem, $D' \cup \{0\}$ is smooth. By the strong maximum principle and embeddedness, $N \cap B_{r/2}(q)$ is smooth. Therefore $N$ is smooth.

**Example 1.** Let $\{\gamma_t\}_{t \geq 0}$ be an immersed curve evolving in $\mathbf{R}^2$ that develops a cusp singularity as in [A1]. The blowup of the singularity is a double density line. Taking the product with a line, the blowup (in the Huisken scaling) is a double density plane.

In this example, $Q = \emptyset$. The example is essentially immersed, because even in the rescaled surfaces $M^\lambda_{-1}$, the crossing approaches the origin as $\lambda \to 0$.

**Example 2.** For $n = 4$, we construct an embedded flow in two pieces such that the blowup is two intersecting planes. Let $\{\gamma^i_t\}_{t \in [0,T]}$ be curves flowing by mean curvature in $\mathbf{R}^3$ that are disjoint except at $t = T$, at which time $\gamma^1_T \cap \gamma^2_T = \{0\}$. Define $M_t = (\gamma^1_t \cup \gamma^2_t) \times \mathbf{R} \subseteq \mathbf{R}^4$. The blowup at $(0, T)$ is the union of two planes:

$$N = (T_0 \gamma^1_T \times \mathbf{R}) \cup (T_0 \gamma^2_T \times \mathbf{R}).$$

Generically, the two planes intersect in a line, but if one curve is a helix and the other some skew line, it can be arranged that the two planes are equal (and $\nu = 2\mathcal{H}^2 \lfloor N$).

An embedded manifold cannot blow up smoothly to two transverse planes, but we have the following conjectured example, inspired by the slow blowup mechanism of Velazquez [V].

**Example 3.** We propose a connected example with $N = \mathbf{C}_z \cup \mathbf{C}_w$ (orthogonal complex planes) as follows. Define $M_t$ such that

$$M_t \cap B_{\varepsilon(t)} \approx \{zw = \varepsilon(t)^2\}, \qquad M_t \setminus B_{\varepsilon(t)} = M'(t) \cup M''(t), \qquad t < T,$$

where $\varepsilon(t) = o(\sqrt{T-t})$. Now $\{zw = \varepsilon(t)^2\}$ is a complex surface, therefore a minimal surface, so $M_0 \cap B_{\varepsilon(t)}$ is moving very slowly. Also, $M'(t)$ (resp. $M''(t)$) is $\mathbf{C}_z$ (resp. $\mathbf{C}_w$) perturbed by a decaying eigenfunction of the linearization of (5) over $\mathbf{C}$. Then $(T-t)^{-1/2} \cdot M_t$ (rescaled coordinates) decays toward $\mathbf{C}_z \cup \mathbf{C}_w$ as $t \to T$.

**Example 4.** We conjecture an example similar to Example 3 but with $\nu = 2\mathcal{H}^2 \lfloor \mathbf{C}_w$ by taking

$$M_t \cap B_{\varepsilon(t)} \approx \{\varepsilon(t)w = z^2\}, \qquad M_t \setminus B_{\varepsilon(t)} \approx M'(t),$$

where $M'(t)$ is a double-valued graph over $\mathbf{C}_w \setminus B_{\varepsilon(t)}$, and is obtained by perturbing $\mathbf{C}_w$ by a double-valued, decaying eigenfunction of the linearization of (5) over $\mathbf{C}_w$.

## §5 A Small $\varepsilon$ Theorem

The following theorem shows that if $\iint |\vec{A}|^2$ is small on parabolic cylinders at all scales, then we have regularity. It generalizes Ecker [E, Th. 2.1] somewhat. Recall that $B_r$ always denotes an open ball.

**Theorem 14 (Small $\varepsilon$ Regularity Theorem).** There is $\varepsilon_5 = \varepsilon_5(n, k)$ such that if $\{M_t^k\}_{t \in [0,1)}$ is a mean curvature flow smoothly immersed in $B_1 \subseteq \mathbf{R}^n$, and

$$(26) \qquad r^{-k} \int_{t-r^2}^{t} \int_{M_t \cap B_r(x)} |\vec{A}|^2 \, d\mu_t \, dt \le \varepsilon^2 \le \varepsilon_5^2$$

whenever $B_r(x) \times [t - r^2, t) \subseteq B_1 \times [0, 1)$, then $M_t \cap B_1$ can be smoothly extended to $B_1 \times \{1\}$ and

(27) $$|\vec{A}(x,t)| \leq C(n,k) \max\left(\frac{1}{1-|x|}, \frac{1}{t^{1/2}}\right) \varepsilon$$

for all $(x, t) \in B_1 \times (0, 1]$.

First, let us express (1) when $M_t$ is a graph. Let $T$ be a $k$-plane in $\mathbf{R}^n$. Let $i, j, \ldots$ be coordinates on $T$, $I, J, \ldots$ coordinates on $T^\perp$, and $\alpha, \beta, \ldots$ coordinates on $R^n$. Let $M_t$ be the graph of $u = (u^I) : T \to T^\perp$, and the image of $U^\alpha = (x^i, u^I(x))$. A tedious calculation reveals that the second fundamental form of $M_t$ pulled back to $T$ satisfies

(32) $$A^\alpha_{ij} \equiv D_i U^\beta D_j U^\gamma A^I_{\beta\gamma} = Q^\alpha_I D^2_{ij} u^I$$

where $Q = (Q^\alpha_I)$ is the orthogonal projection onto $(T_x M_t)^\perp$, restricted to $T^\perp$. Therefore

$$H^\alpha = g^{ij} Q^\alpha_I D_{ij} u^I$$

where $g^{ij} = g^{ij}(Du)$ is the inverse of the pulled-back metric

$$g_{ij} \equiv D_i U^\alpha D_j U^\alpha = \delta_{ij} + D_i u^I D_j u^I.$$

Equation (1) implies the relation $H^\alpha = Q^\alpha_I \partial u^I / \partial t$, since $(\partial u^I/\partial t)$ and $(H^\alpha)$ produce the same motion in the direction normal to $M_t$. Therefore, (1) reads

(28) $$\frac{\partial}{\partial t} u^I = g^{ij}(Du) D^2_{ij} u^I$$

for the graph $u$.

Also

$$|\vec{A}|^2 = g^{ik} g^{jl} Q^\alpha_I D_{ij} u^I Q^\alpha_J D_{kl} u^J.$$

Evidently

$$(\delta^{ij}) \leq (g^{ij}) \leq C(1+|Du|^2)(\delta^{ij}), \qquad c(1+|Du|^2)^{-1}(\delta_{IJ}) \leq (Q_I^\alpha Q_J^\alpha) \leq (\delta_{IJ})$$

in the usual way of ordering matrices, where $\delta$ denotes various identity matrices. Therefore $|Q^{-1}| \leq C(1+|Du|)$ and

(29) $$|\vec{A}| \leq |D^2 u| \leq C(1+|Du|^3)|\vec{A}|$$

where $D^2 u$, $Du$ are measured with respect to the standard flat metric on $T$ and $|\vec{A}|$ is measured with respect to the induced metric $g_{ij}$. Finally

$$c(1+|Du|^{k'})^{-1} U^*(\mathcal{H}^k \lfloor M_t) \leq \mathcal{H}^k \lfloor T \leq U^*(\mathcal{H}^k \lfloor M_t),$$

$k' = \min(k, n-k)$. For convenience, let us write $P_r(x,t) = B_r(x) \times (t-r^2, t]$, $P_r = P_r(0,1)$, $M = \cup M_t \times \{t\}$.

**Lemma 15.** If $\{M_t\}_{t \in [0,1)}$ is a smooth solution of (1) in $B_1$, $0 \in M_1$, and

$$\sup_{M \cap P_1} |\vec{A}(x,t)| \leq 1$$

then

$$|\vec{A}(0,1)|^2 \leq C(n,k) \int_0^1 \int_{M_t \cap B_1} |\vec{A}|^2 \, d\mathcal{H}^2 \, dt.$$

In the following proof let $C$ denote any constant depending on $k$ and $n$.

**Proof.** 1. First we express $M$ as a graph near $(0,1)$. Let $T = T_0 M_1 \subseteq \mathbf{R}^n$, $\Omega_\beta = (T \cap B_\beta) \times [1-\beta^2, 1]$, $C_\beta = (T^\perp)^{-1}(\Omega_\beta)$ (the cylinder over $\Omega_\beta$). Let $M_\beta$ be the component of $M \cap C_\beta$ containing $(0,0)$. Assume for the moment that $M_\beta$ is the graph of a vector function $u : \Omega_\beta \to T^\perp$ (that is, $M_\beta = \{(x + u(x,t), t) : (x,t) \in \Omega_\beta\}$) such that

(30) $$|Du| \leq 1.$$

Then by (29), $|D^2 u| \leq C$. Integrating in $x$,

$$|Du - (Du)| \leq C\beta, \qquad |u - (u) - x \cdot (Du)| \leq C\beta^2,$$

where $(Du)$ denotes the spatial average at a fixed time

$$(Du) = \frac{1}{\omega_k \beta^k} \int_{T \cap B_\beta \times \{t\}} Du \, d\mathcal{H}^k$$

and similarly for $(u)$. Also,

(31) $$\left|\frac{\partial}{\partial t} u\right| \leq C|\vec{H}|(1 + |Du|) \leq C.$$

Since $|Du(0,1)| = |u(0,1)| = 0$, we obtain

$$|Du(\cdot, 1)| \leq C\beta, \qquad |u(\cdot, 1)| \leq C\beta^2,$$

and thus integrating in $t$,

$$|u| \leq C\beta^2, \qquad |x \cdot (Du)| \leq |u| + |(u)| + C\beta^2 \leq C\beta^2,$$

so

$$|Du| \leq |(Du)| + C\beta \leq C\beta.$$

Since $M_\beta$ is a graph satisfying (30) for sufficiently small $\beta$, a continuation argument shows that $M_{\beta_0}$ is a graph over $\Omega_{\beta_0}$ of a function $u$ satisfying

$$|u| + |Du| + |D^2 u| \leq C(n, k)$$

provided $\beta_0 = 1/C$.

2. Differentiating (28) twice and defining $B_{pq}^I = D_{pq}^2 u^I$, we obtain a system of the form

$$\frac{\partial}{\partial t} B_{pq}^I = D_i(a_{ij}(x) D_j B_{pq}^I) + b_i(x) D_i B_{pq}^I + c_{ipqrs}(x) D_i B_{rs}^I$$
$$+ d_{iJ}^I(x) D_i B_{pq}^J + e_{pqrs}(x) B_{rs}^I$$

where, using the boundedness of $|Du| + |D^2u|$, $a$ is uniformly elliptic and $b$, $c$, $d$, $e$ are all bounded with constants depending only on $n$ and $k$. Then Moser iteration [LSU, VII Theorem 2.2] applies to yield

$$\sup |D^2 u| \leq C(n,k) \int_{\Omega_\beta} |D^2u|^2 d\mathcal{H}^k dt.$$

which becomes by (29)

$$|\vec{A}(0,1)| \leq C(n,k) \int_{1-\beta_0^2}^1 \int_{B_{\beta_0}} |\vec{A}|^2 \, d\mathcal{H}^k \, dt,$$

the desired conclusion.

**Proof of Theorem 14.** The proof is a straightforward generalization of [E, Appendix]. By using the Arzela-Ascoli Theorem, it suffices to prove the estimate (27). By approximating with a slightly smaller domain, we may assume that $\{M_t \cap \bar{B}_1\}_{t \in [0,1]}$ is smooth.

1. Define

(33)
$$\Lambda = \sup_{M \cap \bar{P}_1} |\vec{A}(x,t)| \min(1-|x|, t^{1/2})$$

and choose $(x_0, t_0) \in \bar{P}_1$ such that

(34)
$$|\vec{A}(x_0, t_0)| \min(1-|x_0|, t_0^{1/2}) = \Lambda$$

Note that $(x_0, t_0) \in P_1$ since $|\vec{A}|$ is bounded.

If $\Lambda > 1$, define

(35)
$$r \equiv \frac{1}{2|\vec{A}(x_0, t_0)|} \leq \frac{\Lambda}{2|\vec{A}(x_0, t_0)|} = \frac{\min(1-|x_0|, t_0^{1/2})}{2}.$$

Then for $(x,t) \in P_r(x_0, t_0)$ we have

(36)
$$1 - |x| \geq \frac{1-|x_0|}{2}, \qquad t \geq \frac{t_0^2}{4},$$

so that $P_r(x_0, t_0) \subseteq P_1$ and by (33), (34) and (35), (36),

$$|\vec{A}(x,t)| \leq \frac{\min(1-|x_0|, t_0^{1/2})}{\min(1-|x|, t^{1/2})} |\vec{A}(x_0, t_0)| \leq 2|\vec{A}(x_0, t_0)| = \frac{1}{r}$$

for $(x,t) \in P_r(x_0, t_0)$. Then by (35) and Lemma 15 scaled to $P_r(x_0, t_0)$,

$$\begin{aligned} \frac{1}{4} &= r^2 |\vec{A}(x_0, t_0)|^2 \\ &\leq C(n,k) r^{-k} \int_{t-r^2}^{t} \int_{M_t \cap B_r(x)} |\vec{A}|^2 \, d\mathcal{H}^k \, dt \\ &\leq C(n,k) \varepsilon^2 \end{aligned}$$

by hypothesis, which is a contradiction provided we define $\varepsilon_5^2 = 1/8C(n,k)$.

2. Therefore $\Lambda \leq 1$. This time, define

$$r \equiv \frac{\Lambda}{2|\vec{A}(x_0, t_0)|} = \frac{\min(1-|x_0|, t_0^{1/2})}{2}$$

(less rescaling) and as before

$$|\vec{A}(x,t)| \leq \frac{1}{r}$$

for $(x,t) \in P_r(x_0, t_0)$, permitting us to apply Lemma 15 to obtain

$$\begin{aligned} \frac{1}{4} \min(1-|x_0|, t_0^{1/2})^2 |\vec{A}(x_0, t_0)|^2 &= r^2 |\vec{A}(x_0, t_0)|^2 \\ &\leq C(n,k) r^{-k} \int_{t-r^2}^{t} \int_{M_t \cap B_r(x)} |\vec{A}|^2 \, d\mathcal{H}^k \, dt \\ &\leq C(n,k) \varepsilon^2, \end{aligned}$$

so that by (33), (34),

$$\sup_{M \cap \bar{P}_1} |\vec{A}(x,t)| \min(1-|x|, t^{1/2}) \leq 2C(n,k)\varepsilon,$$

as required.